# User Manual for Model-based Imaging Inverse Problem

**Xiaodong Wang**

**wangxiaodong@westlake.edu.cn**

This user manual is intended to provide a detailed description on model-based optimization for imaging inverse problem. These problems can be particularly complex and challenging, especially for individuals without prior exposure to convex optimization or inverse problem theory, like myself. In light of this, I am writing this manual to clarify and systematically organize the mathematical rationale underlying imaging inverse problems. This manual might not be accurate in mathmatical notion but more focus on the logical thinking on how to solve and proceed to solve the problems. If you want to think deep about something, try to raise questions! This manual is seaprated into four sections, aiming to answer the following four questions:

- What is inverse imaging problem?
- Why optimization is used to solve the inverse imaging problem?
- How to solve the optimization problem?
- How to implement the optimization algorithm in real imaging system?

# 1. Inverse problem for computational imaging

> THIS ANSWERS WHAT IS INVERSE IMAGING PROBLEM

Solving imaging inverse problem is to estimate an unknown image $x \in \mathbb{R}^n$ from its corrupted measurements $y \in \mathbb{R}^m$. The acquisition process of the measurement is refered to as a *forward model*, often represented as

$$y = A(x) + e \qquad (1)$$

The forward operator $A : \mathbb{R}^n \to \mathbb{R}^m$ represents the physics of the imaging system, which varies across different imaging modalities. The vector $e \in \mathbb{R}^m$ represents the measurement noise. Table 1 highlights the most popular forward operator used in the context of computational imaging. It can be seen from the table that most inverse imaging problems are linear imaging system, which can be expressed as $y = Ax$. Therefore hereafter the mathmatical derivation in this manual will focus on linear imaging system.

Table 1: Examples of inverse problems in imaging

| Imaging modality | Forward model | Notes |
|---|---|---|
| Denoising | $y = Ix$ | $I$ is the identity matrix. |
| Deblur | $y = h * x = Hx$ | $h$ is a known blur kernel and $*$ denotes convolution ($H$ is a discrete convolution operator). When $h$ is unknown the reconstruction problem is known as blind deconvolution. [1] |
| Superresolution | $y = SBx$ | $S$ is a subsampling operator (identity matrix with missing rows) and $B$ is a blurring operator cooresponding to convolution with a blur kernel. [2] |
| Inpainting | $y = Sx$ | $S$ is a diagonal matrix where $S_{i,i} = 1$ for the pixels that are sampled and $S_{i,i} = 0$ for the pixels that are not. [3] |
| Magnetic resonance imaging (MRI) | $y = SFMx$ | $S$ is a subsampling operator (identity matrix with missing rows), $F$ is the discrete Fourier transform matrix, and $M$ is a diagonal matrix representing a spatial domain multiplication with the coil sensitivity map (assuming a single coil aquisition with Cartesian sampling in a SENSE framework). [4] |
| Computed tomography (CT) | $y = Rx$ | $R$ is the discrete Radon transform. [5] |
| Snapshot compressive imaging (SCI) | $y = Mx$ | $M$ is the sensing matrix related to the 3D mask. [6] |
| Single pixel imaging (SPI) | $y = Mx$ | $M$ is the sensing matrix related to the 3D mask. [7] |
| Non-line-of-sight Imaging (NLOS) | $y = R_t^{-1} H R_z x$ | The matrix $H$ represents the shift-invariant 3D convolution operation, and the matrices $R_t$ and $R_z$ represent the transformation operations applied to the temporal and spatial dimensions, respectively. [8] |
| Structured illumination microscopy (SIM) | $y_i = SHM_i x_i$ | $S$ is a decimation operator with a downsampling factor of two in each dimension, this is required because SIM aims at doubling the lateral resolution. $M_i$ is a diagonal matrix associated to the $i$th illumination patterns. $H$ is a discrete convolution operator.[9] |
| Optical diffraction tomography (ODT) | $y_i = M_i F x$ | $F$ is the discrete Fourier transform matrix, $M_i$ is a diagonal matrix associated to the $i$th illumination patterns. [10] |
| Phase retrival (PR) | $y = \|Ax\|^2$ | $\|\cdot\|$ denotes the absolute value, the square is taken elementwise, and $A$ is a (potentially complex valued) measurement matrix that depends on the application. The measurement |

| Imaging modality | Forward model | Notes |
|---|---|---|
| | | matrix $A$ is often a variation on a discrete Fourier transform matrix. [11] |
| Phase unwrapping (PU) | $y = x - 2\pi k$ | $k$ denotes the map of wrap counts, indicating the number of times a phase value has been wrapped around by $2\pi$. [12] |

## 2. Bayesian approach

> THIS TURNS INVERSE IMAGING PROBLEM INTO OPTIMIZATION PROBLEM

In practice, to recover $x$ from noisy measurement $y$ is nearly impossible since the inverse problem in eq.(1) is often ill-posed (meaning forward operation is a fat matrix ($n \ll m$) in linear imaging system). Therefore, a naive inverting of this forward model will not work. The prior knowledge about the unknown image is needed to further constrain the feasible solution. The detailed mathmatical understanding can be obtained via a *Bayesian approach*. Given the measurement $y$ and forward operator $A$, the Bayesian framework allows us to predict the distribution of a future data $x$ conditioned on $A$ and $y$. To do this, we first obtain the posterior distribution $p_{x|y}$, which can be computed using *Bayers' theorem*

$$p_{x|y}(x|y) = \frac{p_{y|x}(y|x)p_x(x)}{p_y(y)} \propto \underbrace{p_{y|x}(y|x)}_{Likelihood} \underbrace{p_x(x)}_{Prior} \tag{2}$$

Here we use $\propto$ to denote equality after normalization, $p_{y|x}$ is the likelihood function that characters how well the meaurement matches our predictions. $p_x$ is the prior probability that represents our pre-existing knowledge or belief about the true image. The goal of imaging inverse problem is to find an estimation of $x$ that maximizes the posterior distribution $p_{x|y}(x|y)$, which is also known to be *Maximun A Posteriori (MAP)* estimate

$$\widehat{x}_{MAP} = \arg\max_{x}\{p_{x|y}(x|y)\} \tag{3}$$

$$= \arg\min_{x}\{-\log(p_{y|x}(y|x)) - \log(p_x(x))\} \tag{4}$$

$$= \arg\min_{x}\{\underbrace{\mathbf{D}(\mathbf{x})}_{Data-fidelity} + \underbrace{\mathbf{R}(\mathbf{x})}_{Regularizer}\} \tag{5}$$

By using bayesian estimation, the inverse problem is transformed into an optimization problem, with which two terms are included in the objective function, the Data-fidelity term $\mathbf{D}(\mathbf{x})$ and Regularizer term $\mathbf{R}(\mathbf{x})$. $\mathbf{D}(\mathbf{x})$ controls the data consistency to the measurement. $\mathbf{R}(\mathbf{x})$ as a prior term, constrains the solution to feasible distribution.

We assume the noise term in eq.(1) is *addictive white Gaussian noise (AWGN)* with $e \sim N(0, \sigma^2 I)$, that way, the livelihood can be expressed as

$$p_{y|x}(y|x) = \frac{1}{2\pi\sigma^2 I}\exp(-\frac{1}{2\sigma^2}||A(x) - y||_2^2) \tag{6}$$

where $||\cdot||_2^2$ denotes the $l_2 - \text{norm}$. When substituting this equation into eq.(4), we obtain

$$\widehat{x}_{MAP} = \arg\min_{x}\{-\log(p_{y|x}(y|x)) - \log(p_x(x))\} \tag{7}$$

$$= \arg\min_{x}\{\frac{1}{2\sigma^2}||A(x) - y||_2^2 + \mathbf{R}(\mathbf{x})\} \tag{8}$$

$$= \arg\min_{x}\{\frac{1}{2}||A(x) - y||_2^2 + \tau\mathbf{R}(\mathbf{x})\} \tag{9}$$

where $\tau$ is a regularizer parameter to balance the relative strength between $\mathbf{D}(\mathbf{x})$ and $\mathbf{R}(\mathbf{x})$. Therefore, we formulate the inverse problem into an optimization problem. In linear system, by substituting the forward operator $A(x)$ with $Ax$, we have the regularized least-squares

optimization as

$$\widehat{x}_{MAP}^{LS} = \arg\min_x \{\frac{1}{2}||Ax - y||_2^2 + \tau \mathbf{R}(\mathbf{x})\} \quad (10)$$

Hereafter we use $||\cdot||^2$ to denote the $l_2-\text{norm}$ for simplicity. Notice that the data-fidelity term satisfies the $l_2-\text{norm}$ is based on the assumption of gaussian noise, which can also be extended to other noise types. This optimization formulation will be the core idea in the context of computational imaging bacease most imaging system is a low-frequency linear system. Even if it is not, it can be transformed or approximated as a linear one in some way. We will focus one specific linear model, snapshot comoressive imaging in the example section, to demonstrate the following optimization framework.

# 3. Optimization framework

> THIS EXPLAINS HOW TO SOLVE OPTIMIZATION PROBLEM

We reformulate the regularized optimization in eq.(10) as

$$\widehat{x} = \arg\min_x f(x), \ s.t. \ f(x) = g(x) + r(x) \quad (11)$$

where $g(x) = \frac{1}{2}||Ax - y||^2$ and $r(x) = \tau \mathbf{R}(\mathbf{x})$ are data-fidelity term and regularizer term respectively. Now, to solve the above optimization problem, one may use the gradient decent algorithm to iteratively approach the true image.

$$x^k \leftarrow x^{k-1} - \gamma \nabla f(x) \quad (12)$$

where $\nabla$ and $\gamma$ represent the gradient of a function and weight parameter respectively. To compute the gradient of $f(x)$ requires the objective function to be smooth. Nevertheless, in most cases, a large number of regularizers $r(x)$, including $l_1 \text{ norm}$ and TV, are not non-smooth. *Proximal method* is an commonly used method to solve this non-smooth problem. Specfically, *proximal operator* is defined as

$$\text{prox}_{\tau r}(\mathbf{z}) = \arg\min_{\mathbf{x}} \{\frac{1}{2}||\mathbf{x} - \mathbf{z}||^2 + \tau r(\mathbf{x})\} \quad (13)$$

where $\tau$ is a weighting parameter to control the strength of $r(x)$. The *proximal operator* is potent tool for solving optimization problems because in some cases the solution is unique and can be obtained analytically.[13] There are fast and exact methods to compute the proximity operators for a large class of functions, see also the website http://proximity-operator.net.

We will not dive into the details of *proximal operator*, but bear in mind that solving an optimization problem with non-smoothness in eq.(11) is difficult, but solving eq.(13) is accessible. Another interesting note is that the *proximal operator* is equivalent to the forward operator in image denoising problem (also known as the *Moreau proximal mapping*). In the context of regularized optimization problem, all the derivation below is around the problem of how to transform a non-smooth optimization in eq.(11) into an easy-to-use denoising model. Therefore we believe that the fundamental principle behind most regularized optimization is an image denoising problem.

## 3.1 Classical iterative solvers

We survey the following optimization solvers.

- Proximal Gradient Method (PGM)
- Primal Dual (PD)
- Alternating Direction Method of Multipliers (ADMM)
- Half Quadratic Splitting (HQS)
- Generalized Alternating Projection (GAP)

### 3.1.1 Proximal Gradient Method (PGM)

Proximal Gradient Method (PGM), also known as *iterative shrinkage/thresholding algorithm (ISTA)*, consists in alternating an explicit gradient step on a smooth $g$ and implicit gradient descent step on a potentially nonsmooth $r$. [14]

Before diving into detailed derivation of *ISTA*, we introduce the concept of convex function, $L$-Lipschitz continuous gradient and first-order convexity inequality.

**Definition S.1.** A continuously differentiable function $f$ is convex on $\mathbb{R}^n$ if for all $x \in \mathbb{R}^n$ and $y \in \mathbb{R}^n$

$$f(y) \geq f(x) + \nabla f(x)^T(y-x) \tag{14}$$

**Definition S.2.** A function $g$ has $L$-Lipschitz gradient if there exists $L > 0$ such that

$$\frac{L}{2}||x||^2 - g(x) \text{ is convex}, \quad \forall \mathrm{x} \in \mathbb{R}^\mathrm{n} \tag{15}$$

**Definition S.3.** Let $f$ be convex and differentiable function with $\nabla f$ that is $L$-Lipschitz continuous. Then, for all $x \in \mathbb{R}^n$ and $y \in \mathbb{R}^n$

$$f(y) \leqslant f(x) + \nabla f(x)^T(y-x) + \frac{L}{2}||y-x||^2 \tag{16}$$

The PGM assumes the data-fidelity term $g$ is continuously differentiable and has a Lipschitz continuous gradient with constant $L > 0$. By applying the first-order convexity inequality (Definition B.3) to $g$, we obtain the upper bound for the data-fidelity term $g$

$$g(x^+) \leqslant g(x) + \nabla g(x)^T(x^+ - x) + \frac{L}{2}||x^+ - x||^2 \tag{17}$$

By replacing $L$ with $1/\gamma$ ($0 < \gamma < 1/L$), we obtain the upper bound function for $g$

$$q(x^+, q) := g(x) + \nabla g(x)^T(x^+ - x) + \frac{1}{2\gamma}||x^+ - x||^2 \tag{18}$$

Therefore, the optimization problem can be solved by minimizing the sum of upper bound of $g$ and regularizer $r$ at iteration $k \geq 1$

$$x^k = \arg\min_x \{q(x, x^{k-1}) + r(x)\} \tag{19}$$

$$= \arg\min_x \{ g(x^{k-1}) + \nabla g(x^{k-1})^T(x - x^{k-1}) + \frac{1}{2\gamma}||x - x^{k-1}||^2 + r(x)\} \tag{20}$$

$$= \arg\min_x \{\frac{1}{2}||x - (x^{k-1} - \gamma\nabla g(x^{k-1}))||^2 + \gamma r(x)\} \tag{21}$$

$$= \mathrm{prox}_{\gamma \mathrm{r}}(\mathrm{x}^{\mathrm{k}-1} - \gamma\nabla \mathrm{g}(\mathrm{x}^{\mathrm{k}-1})) \tag{22}$$

We now transform the non-smooth optimization into the calculation of proximal operator. By splitting the updates into two steps, we obtain

$$\begin{cases} z^k & \leftarrow x^{k-1} - \gamma\nabla g(x^{k-1}) \\ x^k & \leftarrow \mathrm{prox}_{\gamma \mathrm{r}}(\mathrm{z}^\mathrm{k}) \end{cases} \tag{23}$$

where $\gamma > 0$ is the step size parameter. The above formulation is widely adopted in computational imaging optimization. For linear imaging system, $g(x) = \frac{1}{2}||Ax - y||^2$, we can compute the gradient of $g(x)$ as

$$\nabla g(x) = A^T(Ax - y) \tag{24}$$

Therefore, ISTA solve the linear inverse problem iteratively

$$\begin{cases} z^k & = x^{k-1} - \gamma A^T(Ax^{k-1} - y) \\ x^k & = \mathrm{prox}_{\gamma \mathrm{r}}(\mathrm{z}^\mathrm{k}) \end{cases} \tag{25}$$

where $T$ denotes the transpose operation for real value. PGM can be further extended to an accelerated version by introducing a sequence $q_k$, which is known as *fast iterative shrinkage/thresholding algorithm (FISTA)*

$$\begin{cases} z^k = s^{k-1} - \gamma A^T(Ax^{k-1} - y) \\ x^k = \text{prox}_{\gamma r}(z^k) \\ s^k = x^k + \frac{q_{k-1}-1}{q_k}(x^k - x^{k-1}) \end{cases} \quad (26)$$

where the values for $\{q_k\} = 1$ and $\{q_k\} = \frac{1}{1+\sqrt{1+4q_{k-1}^2}}$ serve as a switch between *ISTA* and *FISTA*. There are various types of variants of ISTA algorithm, such as *two-step iterative shrinkage/thresholding algorithm (TwIST)* [15] and *approximate message passing (AMP)* [16] [17]. On a more general form, *TwIST* iterates by

$$x^k = (1-\alpha)x^{k-2} + (\alpha-\beta)x^{k-1} + \beta \text{ prox}_{\gamma r}(x^{k-1} - \gamma A^T(Ax^{k-1} - y)) \quad (27)$$

where $\alpha > 0$ and $\beta > 0$. The iterative updating can be obtained by splitting the variables as *ISTA* does. It is interesting to notice that when $\alpha = \beta = 1$, the above equation is equavalent to *ISTA*. *FISTA* can also be seen as a special case of *TwIST* by setting $\alpha = 1$. *AMP* extends *ISTA* by adding an extra term to the residual known as the *Onsager correction term* in $z$-updating step

$$z^k = s^{k-1} - \gamma A^T(Ax^{k-1} - y - \frac{1}{\delta}z^{k-1}\langle A^T z^{k-1} + x^{k-1}\rangle) \quad (28)$$

where $\delta = n/m$ is a measure of under-determinacy of the inverse problem, $\langle \cdot \rangle$ denotes the average of a vector.

### 3.1.2 Alternating Direction Method of Multipliers (ADMM)

Instead of calculating the gradient of the data-fidelity term, ADMM solves problems by introducing an auxilliary parameter, hence converting unconstrained optimization into constrained optimization, which in the general form is expressed as

$$\min \quad g(x) + r(z) \quad \text{s.t.} \quad Px + Qz = c \quad (29)$$

where $x \in \mathbb{R}^n, z \in \mathbb{R}^m, P \in \mathbb{R}^{p \times n}, Q \in \mathbb{R}^{p \times m}$ and $c \in \mathbb{R}^p$. When $P = I, Q = -I, c = 0$, ADMM can be simplified as

$$\min \quad g(x) + r(z) \quad \text{s.t.} \quad x = z \quad (30)$$

It can be seen that above equation is simply a constrained version to the optimization problem in eq.(11). The only difference is that a new variable $z$ is introduced, thus separating the variable $x$ into two parts, $x$ and $z$. This variable splitting is the main principle that is commonly adopted in classic optimization. The *augmented Lagrangian* method, originally known as the multipliers method (MM) [18], combines the Lagrangian function and a quadratic penalty term. It is applied to solve a constrained optimization problem iteratively. To solve the constrained optimization problem, ADMM algorithm employs the *augmented Lagrangian* method, which is formulated as

$$L_r(x, z, \mu) = g(x) + r(z) + \mu^T(x-z) + \frac{1}{2\gamma}||x-z||^2 \quad (31)$$

$$= g(x) + r(z) + \frac{1}{2\gamma}||x-z+\gamma\mu||^2 - \frac{1}{2\gamma}||\mu||^2 \quad (32)$$

where $\gamma$ is a regularizer parameter and $\mu \in \mathbb{R}^n$ is the dual variable. By introducing the scaled dual variale $s := \gamma\mu$, we obtain the *scaled augmented Langrangian*

$$L_r(x, z, \mu) = g(x) + r(z) + \frac{1}{2\gamma}||x-z+s||^2 - \frac{1}{2\gamma}||s||^2 \quad (33)$$

ADMM can then be formulated as the joint optimization of three variables, $x, z$ and $s$, which is separated into three subproblems, optimized in an alternating or sequential fashion. This is given by

$$\begin{cases} x^k \leftarrow \arg\min_x L_r(x, z^{k-1}, s^{k-1}) \\ z^k \leftarrow \arg\min_x L_r(x^k, z, s^{k-1}) \\ s^k \leftarrow s^{k-1} + (x^k - z^k) \end{cases} \quad (34)$$

Now let's consider the minimization problem of variable $x$ and $z$ by incorporating proximal operator eq.(13),

$$\arg\min_x L_r(x, z^{k-1}, s^{k-1}) = \arg\min_x g(x) + \frac{1}{2\gamma}||x - z^{k-1} + s^{k-1}||^2 \tag{35}$$

$$= \text{prox}_{\gamma g}(z^{k-1} - s^{k-1}) \tag{36}$$

$$\arg\min_z L_r(x^k, z, s^{k-1}) = \arg\min_z r(z) + \frac{1}{2\gamma}||x^k - z + s^{k-1}||^2 \tag{37}$$

$$= \text{prox}_{\gamma r}(\mathbf{x}^k + \mathbf{s}^{k-1}) \tag{38}$$

Different from PGM that computes the gradient of data-fidelity term $g(x)$, see equation eq.(25), ADMM computes the proximal $\text{prox}_{\gamma g}$ on $g(x)$. For linear imaging system, $g(x) = \frac{1}{2}||Ax - y||^2$, $\text{prox}_{\gamma g}$ can be derived in a closed-form solution by setting the gradient to 0

$$0 = \frac{\partial(\frac{1}{2}||x - z||^2 + \frac{\gamma}{2}||Ax - y||^2)}{\partial x}$$
$$\Rightarrow 0 = x - z + \gamma A^T(Ax - y)$$
$$\Rightarrow z + \gamma A^T y = (I + \gamma A^T A)x$$
$$\Rightarrow x = [I + \gamma A^T A]^{-1}(z + \gamma A^T y) \tag{39}$$

where $I$ is the identity matrix. Therefore, ADMM solve the linear inverse problem iteratively

$$\begin{cases} x^k = [I + \gamma A^T A]^{-1}(z^{k-1} - s^{k-1} + \gamma A^T y) \\ z^k = \text{prox}_{\gamma r}(\mathbf{x}^k + \mathbf{s}^{k-1}) \\ s^k = s^{k-1} + (x^k - z^k) \end{cases} \tag{40}$$

where $\text{prox}_{\gamma r}$ can be seen as a denoising process. Therefore, various denosing algorithms can be used. To conclude, the main takeoff meassage or core idea from ADMM algorithm is that when gradient is not accessible due to the differentiability of the data-fidelity term, variable splitting technique can be employed to formulate two proximal operators, which makes the regularized optimization problem much easier to compute.

### 3.1.3 Primal Dual (PD)

Primal-dual method aims at solving the primal optimization problem[19]:

$$\min_{x \in \mathbb{R}^n} g(x) + r(Dx) \quad \textbf{(primal)} \tag{41}$$

The prior involves a linear operator $D \in \mathbb{R}^{p \times n}$. The modelization of the prior as the composition of a function with a linear operator allows us to model most of the standard penalization such as the anisotropic or isotropic Total Variation (TV). We first give the definition of conjugate function.

**Definition S.3.** The conjugate of a function $f$ is the function $f^*$ defined as

$$f^*(\mu) = \sup_{x \in \mathbb{R}^n}(x^T \mu - f(x)) \tag{42}$$

The corresponding dual formulation of the primal problem is

$$\max_{z \in \mathbb{R}^k} -(g^*(-D^T z) + r^*(z)) \quad \textbf{(dual)} \tag{43}$$

where $*$ denotes the convex conjugate: $r^*(z) = \langle Dx, z \rangle - r(Dx)$, $\langle \rangle$ denotes inner product. As compared to primal problem, the dual problem may be easier to solve, especially when $k$ is much smaller than $n$. The proposed algorithms are $primal - dual$, in that they solve both the primal and the dual problems jointly. The minimation problem is to combine them into the search of a saddle point of the Lagrangian

$$\min_{x \in \mathbb{R}^n} \max_{z \in \mathbb{R}^k} \langle Dx, z \rangle + g(x) - r^*(z) \tag{44}$$

The *primal-dual* is employed in place of the *primal problem* is because *proximal operator* for $r(Dx)$ is not trivial, but we can get *proximal operator* $r^*(x)$ and $g(x)$ easily. A saddle point $(x, z)$ of this min-max function should satisfy

$$\begin{cases} 0 \in D\hat{x} - \partial r^*(\hat{z}) \\ 0 \in D^T\hat{z} + \partial g^*(\hat{x}) \end{cases} \tag{45}$$

Chambolle and Pock (2011) proposed the following algorithm to solve the above problem, which is called *Chambolle-Pock (CP)* or *Primal-Dual Hybrid Gradient (PDHG)*.[20] For stepsizes $\tau > 0$, $\sigma > 0$, and relaxation parameter $\beta \in [0, 1]$, it writes

$$\begin{cases} \hat{x}^k &= \text{prox}_{\tau g}(x^{k-1} - \tau D^T z^{k-1}) \quad \textbf{(primal proximal)} \\ \hat{z}^k &= \text{prox}_{\sigma r^*}(z^{k-1} + \sigma D(2\hat{x}^k - x^{k-1})) \quad \textbf{(dual proximal)} \\ (x^k, z^k) &= (\hat{x}^k, \hat{z}^k) + \beta(x^{k-1}, z^{k-1}) \end{cases} \tag{46}$$

The computation of $\hat{z}^k$ uses the over-relaxed version $2\hat{x}^k - x^{k-1}$ of $x^k$. For linear imaging system, $g(x) = \frac{1}{2}||Ax - y||^2$, the *primal proximal* is esy to be calculated (similar to ADMM in (45) and (46)). It writes

$$\begin{cases} \hat{x}^k &= x^{k-1} - \tau A^T(Ax^{k-1} - y) - \tau D^T z^{k-1} \quad \textbf{(primal proximal)} \\ \hat{z}^k &= \text{prox}_{\sigma r^*}(z^{k-1} + \sigma D(2\hat{x}^k - x^{k-1})) \quad \textbf{(dual proximal)} \\ (x^k, z^k) &= (\hat{x}^k, \hat{z}^k) + \beta(x^{k-1}, z^{k-1}) \end{cases} \tag{47}$$

The core ingredient of *proximal-dual splitting* algorithms is the proximal operator in *dual proximal*. The *Moreau identity* allows us to provide a relation between a function $r$ and its conjugate $r^*$: $\text{prox}_{\sigma r^*}(x) = x - \sigma \text{prox}_{r/\sigma}(x/\sigma)$. Note that we can also formulate the problem in eq.(41) as a *primal problem*: $\min_{x \in \mathbb{R}^n} g(Ax) + r(x)$, which shares the same derivation as described. Primal-dual method is equivalent to ADMM when solving the problem: $\min_{x \in \mathbb{R}^n} g(Ax) + r(x)$ subject to $Ax - y = 0$. Check it out if interested.

### 3.1.4 Half Quadratic Splitting (HQS)

The *Half Quadratic Splitting (HQS)* algorithm takes a similar variable splitting form with ADMM (eq.(30)), the problem of minimizing is equivalent to an constrained optimization problem by introducing an auxiliary variable $z$,

$$\min_{x,z \in \mathbb{R}^n} g(x) + r(z) \quad s.t. \quad x = z \tag{48}$$

Different from ADMM that employs the *augmented Lagrangian* to turn a constrained problem into an unconstrained problem, the method of HQS seeks to minizmize the following cost function by using only quadratic penalty term:

$$\min_{x,z \in \mathbb{R}^n} = g(x) + r(z) + \frac{1}{2\gamma}||x - z||^2 \tag{49}$$

where $\gamma$ is a penalty parameter. This can be solved in an iterative strategy, HQS optimizes $x, z$ in an alternating fashion by solving the following two subproblems separately

$$x^k = \arg\min_x g(x) + \frac{1}{2\gamma}||x - z^{k-1}||^2 \tag{50}$$

$$= \text{prox}_{\gamma g}(z^{k-1}) \tag{51}$$

$$z^k = \arg\min_z r(z) + \frac{1}{2\gamma}||x^k - z||^2 \tag{52}$$

$$= \text{prox}_{\gamma r}(x^k) \tag{53}$$

Again, for linear imaging system, $g(x) = \frac{1}{2}||Ax - y||^2$, we can thus obtain a closed form solution similar to eq.(39)

$$x^k = [I + \gamma A^T A]^{-1}(z^{k-1} + \gamma A^T y) \tag{54}$$

Therefore, HQS solves the linear inverse problem iteratively

$$\begin{cases} x^k &= [I + \gamma A^T A]^{-1}(z^{k-1} + \gamma A^T y) \\ z^k &= \text{prox}_{\gamma r}(x^k) \end{cases} \tag{55}$$

### 3.1.5 Generalized Alternating Projection (GAP)

Generalized alternating projection (GAP)[21] can be recognized as a special case of ADMM. Specfically, GAP introduces an auxiliary parameter $\theta$ and remodels the regularized optimization in eq.(11) as an constrained problem

$$\hat{x}, \hat{\theta} = \arg\min_{x, \theta} \frac{1}{2}||x - \theta||^2 + \gamma r(\theta) \quad s.t. \quad A(x) = y \tag{56}$$

where $\gamma$ is the regularizer parameter. This equation can be solved by alternating updating $\theta$ and $x$

$$x^k \leftarrow \arg\min_x \frac{1}{2}||x - \theta^{k-1}||^2 + \gamma r(\theta^{k-1}) \quad s.t. \quad A(x) = y \tag{57}$$

$$\theta^k \leftarrow \arg\min_\theta \frac{1}{2}||x^k - \theta||^2 + \gamma r(\theta) \tag{58}$$

Now let's consider the minimization problem of variable $x$ and $\theta$. For linear imaging system, $y = Ax$. Given $\theta$, the update of $x$ is simply a *Euclidean projection* of $\theta$ on the linear manifold, this problem in (44) can be modeled as

$$x^k = \arg\min_x \frac{1}{2}||x - \theta^{k-1}||^2 \quad s.t. \quad Ax = y \tag{59}$$

This is equivalent to find the shortest distance between $x$ and $\theta^{k-1}$ on the linear manifold $y = Ax$, which is also a constrained optimization problem. We could get the *Lagrangian function* by introducing a parameter $\mu$,

$$L(x, \mu) = \frac{1}{2}||x - \theta^{k-1}||^2 + \mu(Ax - y) \tag{60}$$

To minimize the *Lagrangian function*, we set the gradient of variables to zero

$$\frac{\partial L(x, \mu)}{\partial x} = x - \theta^{k-1} + \mu A^T = 0 \tag{61}$$

$$\frac{\partial L(x, \mu)}{\partial \mu} = Ax - y = 0 \tag{62}$$

Combining (61) and (62), we obtain

$$y = Ax = A(\theta^{k-1} - \mu A^T) = A\theta^{k-1} - \mu A A^T$$
$$\Rightarrow \mu = -(AA^T)^{-1}(y - A\theta^{k-1}) \tag{63}$$

According to (61), $x = \theta^{k-1} - \mu A^T$. By substituding (63) into $x$, we obtain $x^k$ as

$$x^k = \theta^{k-1} + A^T(AA^T)^{-1}(y - A\theta^{k-1}) \tag{64}$$

Given $x$, the update of $\theta$ is given by

$$\theta^k = \arg\min_\theta \frac{1}{2}||x^k - \theta||^2 + \gamma r(\theta) \tag{65}$$

$$= \text{prox}_{\gamma r}(\mathbf{x}^k) \tag{66}$$

Therefore, GAP solve the linear inverse problem iteratively

$$\begin{cases} x^k = \theta^{k-1} + A^T(AA^T)^{-1}(y - A\theta^{k-1}) \\ \theta^k = \text{prox}_{\gamma r}(\mathbf{x}^k) \end{cases} \tag{67}$$

The linear manifold can also be adaptively adjusted and so (52) can be expressed as an accelerated version of GAP

$$\begin{cases} x^k = \theta^{k-1} + A^T(AA^T)^{-1}(y^{k-1} - A\theta^{k-1}) \\ y^k = y^{k-1} + (y - A\theta^{k-1}) \\ \theta^k = \text{prox}_{\gamma r}(\mathbf{x}^k) \end{cases} \tag{68}$$

Compared with ADMM, GAP and PDM is much simpler since it only uses two steps to update the variables (three steps for accelerated version).

## 3.2 Think of proximal operator as denoiser

As mentioned in **3.1**, the fundamental idea behind the optimization solvers outlined in this article is how to transform the nonsmooth inverse problem into an easy-to-compute *proximal-operator* problem. The *proximal-operator* shares the same mathmatical formulation with an image denoising problem when assuming *AWGN*. This brought a surprising twist in the evolution of regularizers, turning the table and seeking a way to construct a regularization by using denoising algorithms. We introduce three types of denoisers as follows:

- Plug-and-Play (PnP)
- Deep Unfolding Network (DUN)
- Regularization by Denoising (RED)

Based on this observation, researchers replace the proximal operator (or denoiser term, see Eq (25, 40, 47, 55, 67) ) with a neural-network denoiser, which serve as an implicit prior to constrain the feasible solution. Plug-and-Play Priors (PnP) and Deep Unfolding Network (DUN)are two well-known approaches that enable the integration of neural-network-based denoisers within iterative algorithms. In contrast to *PnP* and *DUN* that exploit the implicit prior characterized by a denoiser, *Regularization by Denoising (RED)* uses an explicit image-adaptive Laplacian-based regularization functional, making the overall objective functional clearer and better defined.

### 3.2.1 Plug-and-Play (PnP)

The Plug-and-Play (PnP) methodology, first reported in article [22], proposes to replace the the proximal operator $\text{prox}_{\sigma r}$, within a proximal algorithm with a pretrained image denoiser $\mathbf{D}_\sigma(\cdot)$ ($\sigma$ denotes the denosing strength), such as DnCNN[23] or FastDVDnet[24]. This flexibility enables PnP algorithms to exploit the most effective image denoisers, especially powerful denoising CNNs, leading to their state-of-the-art performance in various imaging tasks. We take PnP-FISTA and PnP-ADMM as examples

**Algorithm 3** PnP-PGM/PnP-APGM
1: **input:** $x^0 = s^0 \in \mathbb{R}^n$, $\gamma > 0$, $\sigma > 0$, and $\{q_k\}_{k \in \mathbb{N}}$
2: **for** $k = 1, 2, \ldots$ **do**
3: $\quad z^k \leftarrow s^{k-1} - \gamma \nabla g(s^{k-1})$
4: $\quad x^k \leftarrow \mathbf{D}_\sigma(z^k)$
5: $\quad s^k \leftarrow x^k + ((q_{k-1} - 1)/q_k)(x^k - x^{k-1})$
6: **end for**

**Algorithm 4** PnP-ADMM
1: **input:** $x^0 \in \mathbb{R}^n$, $s^0 = 0$, $\gamma > 0$, and $\sigma > 0$
2: **for** $k = 1, 2, \ldots$ **do**
3: $\quad z^k \leftarrow \text{prox}_{\gamma g}(x^{k-1} - s^{k-1})$
4: $\quad x^k \leftarrow \mathbf{D}_\sigma(z^k + s^{k-1})$
5: $\quad s^k \leftarrow s^{k-1} + (z^k - x^k)$
6: **end for**

Figure 1. Plug-and-play image restoration based on PGM and ADMM respectively.[25]

### 3.2.2 Deep Unrolling Network (DUN)

Perhaps the most confusing concept to PnP in the model-based deep learning is so-called algorithm unrolling methods, which explicitly unroll/unfold iterative optimization algorithms into learnable deep architectures. In this way, the denoiser prior and other internal parameters are treated as trainable parameters, meanwhile, the number of iterations has to be fixed to enable end-to-end training.

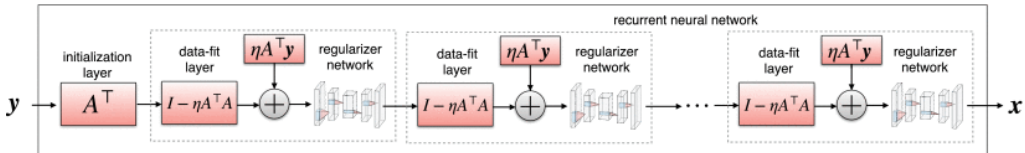

Figure 2. Deep unrolling network based on PGM.[26]

### 3.2.3 Regularization by Denoising (RED)

The RED framework, first proposed in [27], is an alternative scheme where the denoiser $D_\sigma(\cdot)$ can sometimes lead to an explicit regularization function. Specifically, RED suggests to use the following as the regularization function:

$$r(x) = \frac{\lambda}{2} x^T (x - f(x)) \tag{69}$$

where $f(\cdot)$ is a denoiser, $\lambda$ is the regularizer parameter. The advantage of using an explicit denoiser is that various optimization algorithms can be used to solve the following linear inverse problem

$$\hat{x} = \arg\min_x E(x) \tag{70}$$

$$= \arg\min_x g(x) + r(x) \tag{71}$$

$$= \arg\min_x g(x) + \frac{\lambda}{2}x^T(x - f(x)) \tag{72}$$

where $E(x) = g(x) + r(x)$ is the enegy/objective function to be optimized. It is proven numerically in the article [28] that the image denoiser $f(\cdot)$ is *locally homogeneous*, which states

$$(1+\epsilon)f(x) = f((1+\epsilon)x) \quad \forall x \tag{73}$$

for subficiently small $\epsilon \in \mathbb{R}\backslash 0$. Such denoiser obeys the gradient rule

$$\nabla r(x) = \lambda(x - f(x)) \tag{74}$$

If the above equation holds, then any minimizer of the variational objective function requires the gradient to be 0, which writes

$$0 = \nabla g(x) + \lambda(x - f(x)) \tag{75}$$

For linear imaging system, the gradient expression of $E(x)$ is

$$\nabla E_x(x) = A^T(Ax - y) + \lambda(x - f(x)) \tag{76}$$

We now turn to show several options for using this in order to solve a general inverse problem. We list two examples for better illustration:

1. Gradient Descent Method, RED-GD for short

2. ADMM, RED-ADMM for short

1. For Gradient Descent Methods, given the gradient of the energy function $E(x)$, the Steepest-Descent (SD) can be used to update the formula

$$\begin{aligned} x^k &= x^{k-1} - \mu \nabla E_x(x) |_{x^{k-1}} \\ &= x^{k-1} - \mu(A^T(Ax^{k-1} - y) - \lambda(x^{k-1} - f(x^{k-1}))) \end{aligned} \tag{77}$$

where $\mu > 0$ is the step size. Therefore the updating of RED-SGD can be reformulate as

$$\begin{cases} \hat{x}^k &= f_\sigma(x^{k-1}) \\ x^k &= \hat{x}^{k-1} - \mu(A^T(A\hat{x}^{k-1} - y) - \lambda(\hat{x}^{k-1} - f(\hat{x}^{k-1}))) \end{cases} \tag{78}$$

2. For ADMM method, given the gradient of the energy function $E(x)$, the ADMM algorithm introduces an auxiliary parameter $z$ and formulates the constrained optimization problem

$$\arg\min_x \frac{1}{2}\|Ax - y\|^2 + \frac{\lambda}{2}z^T(z - f(z)) \quad \text{s.t.} \quad x = z \tag{79}$$

Following ADMM optimization pipeline, we thus obtaining the iterative updating of $x$ and $z$

$$x^k = \arg\min_x \frac{1}{2}\|Ax - y\|^2 + \frac{1}{2\gamma}\|x - z^{k-1} + s^{k-1}\|^2 \tag{80}$$

$$z^k = \arg\min_z \frac{\lambda}{2}z^T(z - f(z)) + \frac{1}{2\gamma}\|x^k - z + s^{k-1}\|^2 \tag{81}$$

The updating of $x$ can be solved with a closed-form solution $x^k = [I + \gamma A^T A]^{-1}(z^{k-1} - s^{k-1} + \gamma A^T y)$. Hence by setting the gradient of the above expression to zero leads to equation

$$0 = \lambda(z - f(z)) - \frac{1}{\gamma}(x^k - z + s^{k-1}) \tag{82}$$

which can be solved iteratively using the fixed-point strategy, by

$$0 = \lambda(z^k - f(z^{k-1})) - \frac{1}{\gamma}(x^k - z^k + s^{k-1})$$
$$\Rightarrow z^k = \frac{\gamma}{1+\gamma\lambda}(\lambda f(z^{k-1}) + \frac{1}{\gamma}(x^k + s^{k-1})) \tag{83}$$

Therefore we can simplify the above equation and write down the form of iterative updating

$$\begin{cases} x^k &= [I + \gamma A^T A]^{-1}(z^{k-1} - s^{k-1} + \gamma A^T y) \\ \hat{z}^{k-1} &= f_\sigma(z^{k-1}) \\ z^k &= \frac{\gamma}{1+\gamma\lambda}(\lambda f(\hat{z}^{k-1}) + \frac{1}{\gamma}(x^k + s^{k-1})) \\ s^k &= s^{k-1} + (x^k - z^k) \end{cases} \tag{84}$$

# 4. Example: snapshot compressive imaging

> THIS GIVES AN EXAMPLE FOR SOLVING THE INVERSE PROBLEM

## 4.1 Review on classical iterative solvers

Let us review what we have discussed for a linear imaging system. The objective is to solve the following optimization problem

$$\hat{x} = \arg\min_x \frac{1}{2}\|Ax - y\|^2 + r(x)$$

The iterative updating of parameters is summarized below (see Eq (25, 40, 47, 55, 67) ).

(1) Proximal Gradient Method (PGM)

$$\begin{cases} z^k &= x^{k-1} - \gamma A^T(Ax^{k-1} - y) \\ x^k &= \text{prox}_{\gamma r}(z^k) \end{cases}$$

(2) Alternating Direction Method of Multipliers (ADMM)

$$\begin{cases} x^k &= [I + \gamma A^T A]^{-1}(z^{k-1} - s^{k-1} + \gamma A^T y) \\ z^k &= \text{prox}_{\gamma r}(x^k + s^{k-1}) \\ s^k &= s^{k-1} + (x^k - z^k) \end{cases}$$

(3) Primal Dual (PD)

$$\begin{cases} \hat{x}^k &= x^{k-1} - \tau A^T(Ax^{k-1} - y) - \tau D^T z^{k-1} \quad (\textbf{primal proximal}) \\ \hat{z}^k &= \text{prox}_{\sigma r^*}(z^{k-1} + \sigma D(2\hat{x}^k - x^{k-1})) \quad (\textbf{dual proximal}) \\ (x^k, z^k) &= (\hat{x}^k, \hat{z}^k) + \beta(x^{k-1}, z^{k-1}) \end{cases}$$

(4) Half Quadratic Splitting (HQS)

$$\begin{cases} x^k &= [I+\gamma A^T A]^{-1}(z^{k-1}+\gamma A^T y) \\ z^k &= \text{prox}_{\gamma r}(\mathbf{x}^k) \end{cases}$$

(5) Generalized Alternating Projection (GAP)

$$\begin{cases} x^k &= \theta^{k-1} + A^T(AA^T)^{-1}(y - A\theta^{k-1}) \\ \theta^k &= \text{prox}_{\gamma r}(\mathbf{x}^k) \end{cases}$$

Now we know that the *proximal operator* can be replaced by off-the-shelf denoisers, which is represented by implicit neural network or explicit denoising equivalence. The main issue confronting us right now is the quadratic term with relation to the *forward operator* $A$. As we can see in most iterative updating steps, there are mainly three steps involving calculation upon $A$, and they are

- $Ax$ and $A^T y$ in almost all algorithms
- $(AA^T)^{-1}$ in GAP
- $(I+\gamma A^T A)^{-1}$ in ADMM and HQS

In linear imaging system, the *forward operator* is a huge fat matrix, which means more collums than rows. Therefore the arithmetric related to $A$ is very difficult. One direct way to solve the above three arithmetric is to save the matrix $A$ in local space but this could be enoumously inefficient. Fortunately in some specifc cases we do not need to save and do the arithmetric manually but rely on the characteristic of the matrix for efficient computing. We take snapshot compressive imaging (SCI) as an example.

## 4.2 Brief introduction of snapshot compressive imaging

SCI utilizes a 2D detector to capture high-dimensional (≥3D) data. The forward model for SCI system is expressed as

$$y = Ax + e$$

where $y$ is the vectorized measurement, $x = [x_1^T, ..., x_N^T]^T$ is the desired 3D data we are trying to reconstruct, with $x_k$ represents the vectorization of the $k$th frame (by stacking columns). The sensing matrix $\Phi$ is expressed as

$$A = [D_1, \cdots, D_{N_t}] \tag{85}$$

where $D_k = Diag(\text{vec}(M_k))$ is a diagonal matrix with diagonal elements in the vectorized form of modulating mask $M_k$. The explict form of the sensing matrix $A$ is shown in Figure 3.

Figure 3. The forward model (denoted by $\Phi$) of SCI system. Please refer to the original paper for more details.[29]

In an SCI system, the measurement is a 2D image $Y \in \mathbb{R}^{N_x \times N_y}$. Each frame in the 3D signal $X \in \mathbb{R}^{N_x \times N_y \times N_t}$ is modulated by a different mask, resulting in the sensing matrix $A \in \mathbb{R}^{N_x N_y \times N_x N_y N_t}$. Each element in $Y$ is a weighted sum of the corresponding elements in each frame of $X$. Obviously, the compressive sampling rate (CSr) of SCI is $1/N_t$.

## 4.3 Step-by-step coding for snapshot compressive imaging

We start by considering simple algorimetric operations of the forward model $A$. For $Ax$ and $A^T y$, instead of saving the matrix $A$ and do the arithmetic, $Ax$ is calculated by the sum of element-wise multiplication of the mask and scenes (spectrums), and $A^T y$ is calculated by element-wise multiplication of the mask $A$ and the measurement $y$. We showcase the detailed coding on $Ax$ and $A^T y$ with python language.

```python
def H(x, Phi):
    '''
    Forward model of snapshot compressive imaging (SCI), where multiple coded frames are collapsed into a snapshot measurement.

    Parameters
    ----------
    x : 3D ndarray
        Start point (initialized value) for the iteration process of the reconstruction.
    Phi : three-dimensional (3D) ndarray of ints, uints or floats, omitted
        Input sensing matrix of SCI with the third dimension as the time-variant,
        spectral-variant, volume-variant, or angular-variant masks,
        where each mask has the same pixel resolution as the snapshot measurement.
    H : function
        Forward model of SCI, where multiple encoded frames are collapsed into a single measurement.
    '''
    # (nrow, ncol, nmask) = x.shape
    # (nrow, ncol, nmask) = Phi.shape
    return np.sum(x*Phi, axis=tuple(range(2,Phi.ndim)))  # element-wise product

def Ht(y, Phi):
    '''
    Tanspose of the forward model of SCI.

    Parameters
    ----------
    y : two-dimensional (2D) ndarray of ints, uints or float
        Input single measurement of the snapshot compressive imager (SCI).
    Phi : three-dimensional (3D) ndarray of ints, uints or floats, omitted
        Input sensing matrix of SCI with the third dimension as the time-variant,
        spectral-variant, volume-variant, or angular-variant masks,
        where each mask has the same pixel resolution as the snapshot measurement.
    Ht : function
        Tanspose of the forward model of SCI.
    '''
    # (nrow, ncol, 1) = y.shape
    # (nrow, ncol, nmask) = Phi.shape
    return np.multiply(np.repeat(y[:,:,np.newaxis],Phi.shape[2],axis=2), Phi)
```

**We now turn to GAP as an example to illustrate the algorithm implementation.** Recall that in GAP framework, $x^k = \theta^{k-1} + A^T(AA^T)^{-1}(y - A\theta^{k-1})$. It is easy to see that the difficulty lies in the calculation of the matrix inversion on $AA^T$. For $(AA^T)^{-1}$, $AA^T$ is a diagonal matrix and thus the inversion operation of $AA^T$ is also a diagonal matrix. And $AA^T$ is the sum of the square of 3D masks. Here is the derivation:

$$AA^T = \begin{bmatrix} D_1, \cdots, D_{N_t} \end{bmatrix} \begin{bmatrix} D_1^T \\ \vdots \\ D_{N_t}^T \end{bmatrix}$$

$$= D_1 D_1^T + \cdots + D_{N_t} D_{N_t}^T \tag{86}$$

We know that $AA^T$ is a diagonal matrix, with each entry of diagonal value being the square sum of masks. This is expressed as $AA^T = diag(vec(M_1^2 + \cdots + M_{N_t}^2))$. Recall that in GAP framework, $x^k = \theta^{k-1} + A^T(AA^T)^{-1}(y - A\theta^{k-1})$, the iterative updating in GAP algorithm can be implemented as follows

```python
Phi_sum = np.sum(Phi**2,2) # calculation of AA^T
Phi_sum[Phi_sum==0]=1 # make sure non-zero if devided

def gap_denoise(denoiser,theta,y,yy,Phi,Phi_sum,sigma,accelerate=True):
    '''
    GAP framework.

    Parameters
    ----------
    y : Two-dimensional (2D) ndarray of ints, uints or floats
        Input single measurement of the snapshot compressive imager (SCI).
    Phi : three-dimensional (3D) ndarray of ints, uints or floats, omitted
        Input sensing matrix of SCI with the third dimension as the
        time-variant, spectral-variant, volume-variant, or angular-variant
        masks, where each mask has the same pixel resolution as the snapshot
        measurement.
    Phi_sum : 2D ndarray,
        Sum of the sensing matrix `Phi` along the third dimension.
    H : function
        Forward model of SCI, where multiple encoded frames are collapsed into
        a single measurement.
    Ht : function
        Transpose of the forward model.
    theta: 3D ndarray
        Introduced parameter for updating x in accelerated GAP.
    accelerate : boolean, optional
        Enable acceleration in GAP.
    yy: Two-dimensional (2D) ndarray of ints, uints or floats
        Introduced parameter for acceleration in accelerated GAP.
    denoiser : string, optional
        Denoiser used as the regularization imposing on the prior term of the
        reconstruction.
    sigma : One-dimensional (1D) ndarray of ints, uints or floats
        Input noise standard deviation for the denoiser if and only if noise
        estimation is disabled(i.e., noise_estimate==False). The scale of sigma
        is [0, 255] regardless of the the scale of the input measurement and
        masks.

    Returns
    -------
    x : 3D ndarray
        Reconstructed 3D scene captured by the SCI system.

    '''
    yb = H(theta,Phi)
    if accelerate:
        yy = yy + (y-yb)
        x = theta + (Ht((yy-yb)/Phi_sum,Phi))
    else:
        x = theta + (Ht((yy-yb)/Phi_sum,Phi))
    theta = denoiser(x,sigma)
    return x
```

**Another example is ADMM algorithm.** The ADMM algorithm solves the $x^k = (I + \gamma A^T A)^{-1}(z^{k-1} - s^{k-1} + \gamma A^T y)$. For $(I + \gamma A^T A)^{-1}$, since $A$ is a fat matrix, $(I + \gamma A^T A)$ will be a large matrix and thus the matrix inversion formula is employed to simplify the calculation:

$$(I + \gamma A^T A)^{-1} = I - A^T(I + \gamma A A^T)^{-1} A \gamma \qquad (87)$$

This is calculated effciently using the matrix inversion lemma (Woodbury Formula)[30], which gives $(W + UV^T)^{-1} = W^{-1} - [W^{-1}U(I + V^T W^{-1} U)^{-1} V^T W^{-1}]$. We can see from the above equation that another matrix inversion $(I + \gamma A A^T)^{-1}$ is introduced, yet it much easier to be calculated since $AA^T$ can be calculated efficiently, as discussed before. This is achieved by

$$(I + \gamma A A^T)^{-1} = diag\{\frac{1}{\gamma a_1 + 1}, \cdots, \frac{1}{\gamma a_m + 1}\} \qquad (88)$$

where $AA^T \stackrel{\text{def}}{=} diag(a_1, \cdots, a_m)$ ($m = N_t$ in previous SCI forward model). Now if we add the desired 3D data $x$ and consider ADMM algorithm, $x^k = (I + \gamma A^T A)^{-1}(z^{k-1} - s^{k-1} + \gamma A^T y)$. Let $y = [y_1, \cdots, y_m]^T$ and $[A(z^{k-1} - s^{k-1})]_i$ denotes the $i$th element of the vector $A(z^{k-1} - s^{k-1})$, this is further expanded as [31]

$$\begin{aligned}(I + \gamma A^T A)^{-1}(z^{k-1} - s^{k-1} + \gamma A^T y) &= (I - A^T(I + \gamma A A^T)^{-1} A \gamma)(z^{k-1} - s^{k-1} + \gamma A^T y) \\ &= (z^{k-1} - s^{k-1} + \gamma A^T y) - A^T(I + \gamma A A^T)^{-1} A \gamma (z^{k-1} - s^{k-1} + \gamma A^T y) \\ &= z^{k-1} - s^{k-1} + \gamma A^T (y - \underbrace{(I + \gamma A A^T)^{-1}}_{diag\{\frac{1}{\gamma a_1+1},\cdots,\frac{1}{\gamma a_m+1}\}} A(z^{k-1} - s^{k-1}) - \underbrace{\gamma(I + \gamma A A^T)^{-1} A A^T}_{diag\{\frac{\gamma a_1}{\gamma a_1+1},\cdots,\frac{\gamma a_m}{\gamma a_m+1}\}} y) \\ &= z^{k-1} - s^{k-1} + \gamma A^T \left[\frac{y_1 - [A(z^{k-1}-s^{k-1})]_1}{\gamma a_1+1}, \cdots, \frac{y_n - [A(z^{k-1}-s^{k-1})]_m}{\gamma a_m+1}\right]^T \end{aligned} \qquad (89)$$

The above equation can be simplified as

$$x^k = z^{k-1} - s^{k-1} + \gamma A^T(y - A(z^{k-1} - s^{k-1})) \oslash \text{Diag}(\gamma AA^T + I) \qquad (90)$$

where $\text{Diag}(\cdot)$ extracts and vectorizes the diagonal elements of the ensued matrix. $\oslash$ denotes the element-wise division or Hadamard division. Therefore we can write down the iterarive updating codes in ADMM framework

```python
def admm_denoise(denoiser,z,y,s,Phi,Phi_sum,gamma,sigma):
    '''
    ADMM framework.

    Parameters
    ----------
    H : function
        Forward model of SCI, where multiple encoded frames are collapsed into a single measurement.
    Ht : function
        Tanspose of the forward model of SCI.
    Phi : three-dimensional (3D) ndarray of ints, uints or floats, omitted
        Input sensing matrix of SCI with the third dimension as the time-variant,
        spectral-variant, volume-variant, or angular-variant masks,
        where each mask has the same pixel resolution as the snapshot measurement.
    Phi_sum : 2D ndarray
        Sum of the sensing matrix `Phi` along the third dimension.
    y : two-dimensional (2D) ndarray of ints, uints or float
        Input single measurement of the snapshot compressive imager (SCI).
    z : 3D ndarray
        Introduced auxilliary parameter in ADMM.
    s : 3D ndarray
        Introduced scaled dual variale in ADMM.
    denoiser : string, optional
        Denoiser used as the regularization imposing on the prior term of the
        reconstruction.
    gamma : float, optional
        Parameter in the ADMM projection, where more noisy measurements require
        greater gamma.
    sigma : one-dimensional (1D) ndarray of ints, uints or floats
        Input noise standard deviation for the denoiser if and only if
        noise estimation is disabled(i.e., noise_estimate==False).
        The scale of sigma is [0, 255] regardless of the the scale of the input measurement and masks.

    Returns
    -------
    x : 3D ndarray
        Reconstructed 3D scene captured by the SCI system.
    '''

    yb = H(z-s, Phi)
    x = (z-s) + gamma*(Ht((y-yb)/(Phi_sum*gamma+1),Phi))
    z = x-s
    z0 = z
    with torch.no_grad():
        z = denoiser(z, sigma)
    s = s + (x-z)
    return x
```

**End of this Mannual** The prevalent iterative methods employed in solving imaging inverse problems has been surveyed in this manual. While mathematical rigor and detailed derivations may not be the primary focus of this manual, I aim to clarify the overarching framework and knowledge hierarchy inherent to these problems. It is my sincere hope that this guide will serve as a valuable resource on your journey into the realm of imaging inverse problem studies. Should you seek opportunities for academic collaboration or spot areas requiring corrections within this manual, please do not hesitate to reach out.